\documentclass[11pt]{article}
\usepackage{graphics}
\usepackage[demo]{graphicx}
\usepackage{epsfig}
\usepackage{pstricks}
\usepackage[normal]{subfigure}
\usepackage[latin1]{inputenc}
\usepackage[english,francais]{babel}
\usepackage{relsize,exscale}
\usepackage{makeidx}
\usepackage{enumitem}
\usepackage{amsfonts,amssymb,amsmath}
\usepackage{graphicx}
\usepackage{color}
\usepackage{multirow}
\usepackage{mathrsfs}
\usepackage[normalem]{ulem}
\usepackage{cancel}
\usepackage{bbm}
\newenvironment{prooff}{{\it Proof :}}{\hfill\rule{2mm}{2mm}\vskip3mm\par}
\newtheorem{theorem}{Theorem}[section]
\newtheorem{lemma}[theorem]{Lemma}

\newtheorem{e-definition}[theorem]{Definition\rm}
\newtheorem{remark}{\it Remark\/}

\usepackage{color}
\definecolor{dred}{rgb}{0.92,0,0}
\definecolor{dgreen}{rgb}{0,0.92,0}
\definecolor{dblue}{rgb}{0,0,0.92}
\definecolor{dyellow}{rgb}{0.95,0.95,0}

\newcommand{\R}{\mathbb{R}}
\newcommand{\N}{\mathbb{N}}
\def\D{\displaystyle}
\newcommand{\hs}{\hspace{0.1cm}}

\newcommand{\sa}{\\ [0.2cm]}
\usepackage{multirow}
\graphicspath{
{./Figures/}
{./}
}
\title{Generalized Beta Prime Distribution \\ Applied to Finite Element Error Approximation}
\author{Jo\"el Chaskalovic \thanks{D'Alembert,
Sorbonne University, Paris, France, (\emph{email}: jch1826@gmail.com)}
\qquad
Franck Assous
\thanks{
Department of Mathematics, Ariel University, 40700 Ariel, Isra\"el, (\emph{email}: franckassous55@gmail.com).}
}
\date{}
\begin{document}
\maketitle
\selectlanguage{english}
\begin{abstract}
\noindent In this paper we propose a new generation of probability laws based on the generalized Beta prime distribution to estimate the relative accuracy between two Lagrange finite elements $P_{k_1}$ and $P_{k_2}, (k_1<k_2)$. Since the relative finite element accuracy is usually based on the comparison of the asymptotic speed of convergence when the mesh size $h$ goes to zero, this probability laws highlight that there exists, depending on $h$, cases such that $P_{k_1}$ finite element is more likely accurate than the $P_{k_2}$ one.
To confirm this feature, we show and examine on practical examples, the quality of the fit between the statistical frequencies and the corresponding probabilities determined by the probability law. Among others, it validates, when $h$ moves away from zero, that finite element $P_{k_1}$ may produces more precise results than a finite element $P_{k_2}$ since the probability of the event "$P_{k_1}$ \emph{is more accurate than} $P_{k_2}$" consequently increases to become greater than 0.5. In these cases, $P_{k_2}$ finite elements are more likely overqualified.
\end{abstract}
\noindent {\em keywords}: Error estimates, Finite elements, Bramble-Hilbert lemma, Probability.
%
%\vspace{-0.8cm}
%
\section{Introduction}\label{intro}
\noindent Recently, we proposed in \cite{CMAM1}, \cite{CMAM2}, \cite{arXiv_Wmp} and \cite{ArXiv_JCH}  a new vision to consider the error estimate applied to finite elements approximation. Mainly, we derived two probability laws by considering the approximation error as a random variable whose support is determined by the considered error estimate. As a consequence, the relative accuracy between two Lagrange finite elements $P_{k_1}$ and $P_{k_2}, (k_1<k_2),$ is therefore analyzed as random variable too.\sa
This new point of view enabled us to get some distance with the classical results which usually compares the relative accuracy between two finite elements $P_{k_1}$ and $P_{k_2}, (k_1<k_2),$ by the help of the asymptotic speed of convergence when the mesh size $h$ goes to zero to conclude: Since $h^{k_2}$ goes faster to zero than $h^{k_1}$, then the $P_{k_2}$ finite element is more accurate the $P_{k_1}$ one.\sa
The probability laws we derived confirmed our suspicion, (see also our previous different approaches in \cite{AsCh11}, \cite{AsCh13} and \cite{AsCh16}), that when $h$ is set to a \underline{fixed} value, one cannot affirm the same asymptotic conclusion we recalled above. Indeed, since in the error estimates, the upper bound of the approximation error is constituted by an unknown constant which depends, among others, on a given semi-norm of the unknown exact solution, the numerical comparison between the two error approximations associated to $P_{k_1}$ and $P_{k_2}$ cannot be achieved.\sa
As a consequence, to determine the smallest of the two concerned approximation errors is an open-ended question which remains. \sa
From this starting point and by considering a given approximation error as a positive number whose position is unknown within the interval determined by the upper bound of the error estimate, we considered this position as the result of a random variable since the approximation error depends on the approximation, and so, on quantitative uncertainties generated by the process of the mesh generator. \sa
Aside the new insights we got from these probability laws, we also implemented practical cases \cite{Fitting Stat_Heaviside} to appreciate the quality of the fit between these probability laws and the corresponding statistical frequencies. \sa
There, despite the fact we showed that the two probability laws globally behave well like the corresponding statistical frequencies, the fit was not enough precise. Then, we identified the reasons of this unsatisfactory result, basically due to the too much rigidity of the probabilistic assumptions we considered to derive these laws. \sa
It's the reason why we have developed a new generation of probabilistic model based on the generalized Beta prime distribution which enables us to derive, under probabilistic acceptable hypothesis, the probability law of the relative accuracy between two Lagrange finite elements. \sa
In this paper, we will motivate the probabilistic framework we build to get this law and we will show how it fits well with several examples. This is a significant confirmation to assess the relevance of considering the error estimates like random variables in a suitable probabilistic environment. \sa
The paper is organized as follows. We recall in Section \ref{Geo_and_Proba} the mathematical problem we consider and a corollary of Bramble-Hilbert lemma from which we derived the previous probabilistic laws. In Section \ref{First_Comparison}, we show a typical result we got between numerical statistics and these probability laws. Then, in Section \ref{New probabilistic_law}, we derive the new probability law which evaluates the relative error accuracy between two finite elements $P_{k_1}$ and $P_{k_2}, (k_1<k_2).$ by the help of the generalized Beta prime distribution. Finally, in Section \ref{Stat_Proba} we show with several examples the appropriateness of fit between the Generalized Beta prime probabilistic law and the corresponding statistical frequencies. Concluding remarks follow.
\section{The abstract problem and the corresponding finite element error estimate}\label{Geo_and_Proba}
\begin{figure}[htb]
  \centering
  \includegraphics[width=10.cm]{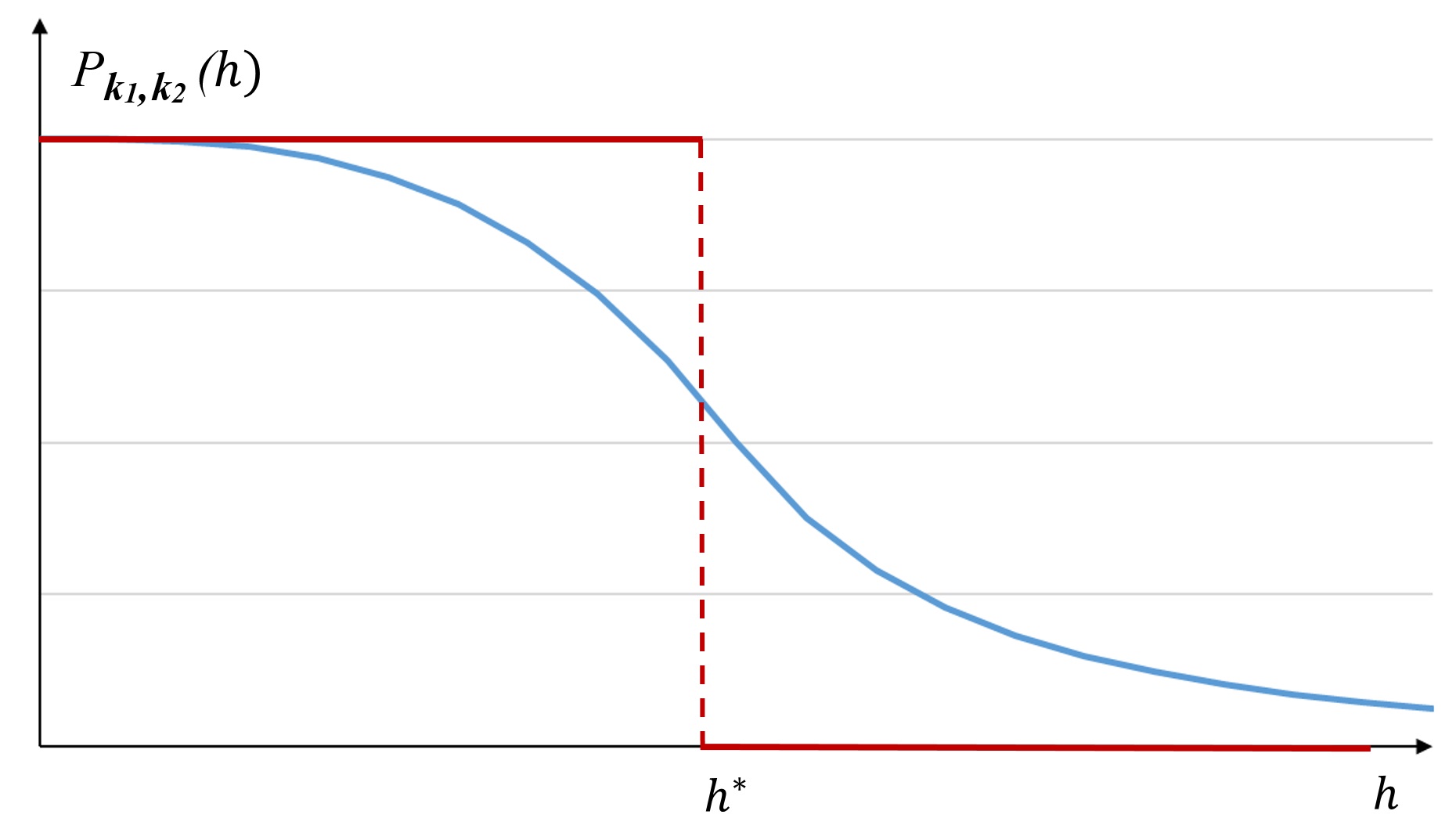}
  \caption{General shape of probabilistic law (\ref{CMAM1})-(\ref{CMAM2}) together with the limit Heaviside case (\ref{Heaviside_Prob}).} \label{Sigmoid_Curve}
\end{figure}

\noindent We consider an open bounded and non empty subset $\Omega$ of $\R^{n}$, and we denote by $\Gamma$ its boundary assumed to be $C^1- $piecewise. We also introduce an Hilbert space $V$ endowed with a norm $ \left\|.\right\|_{V}$, and a bilinear, continuous and $V-$elliptic form $a(\cdot,\cdot)$  defined on $V \times V$. Finally, $l(\cdot)$ denotes a linear continuous form defined on~$V$.\sa
Let $u \in V$ be the unique solution to the second order elliptic variational formulation
\begin{equation}\label{VP}
\left\{
\begin{array}{l}
\mbox{Find } u \in   V \mbox{ solution to:} \\[0.1cm]
a(u,v) = l(v), \quad\forall v \in V\,.
\end{array}
\right.
\end{equation}
In this paper, we will restrict ourselves to the simple case where $V$ is the usual Sobolev space of distributions $H^1(\Omega)$.  More general cases can be found in \cite{ChAs20}.\sa
Let us introduce now $V_h$ a finite-dimensional subset of $V$, and consider $u_{h}\in V_{h}$ an approximation of $u$, solution  to the approximate variational formulation
\begin{equation}\label{VP_h}
\left\{
\begin{array}{l}
\mbox{Find } u_{h} \in   V_h \mbox{ solution to:} \\[0.1cm]
a(u_{h},v_{h}) = l(v_{h}),\quad \forall v_{h} \in V_h.
\end{array}
\right.
\end{equation}
In what follows, we are interested in evaluating error bounds for finite element methods. Hence, we first assume that the domain $\Omega$ is exactly covered by a mesh ${\mathcal T}_h$ composed by $N_{s}$ n-simplexes $K_{j}, (1 \leq j \leq N_{s}),$ which respects classical rules of regular discretization, (see for example \cite{ChaskaPDE} for the bidimensional case or \cite{RaTho82} in $\R^n$). We also denote by $P_k(K_{j})$ the space of polynomial functions defined on a given n-simplex $K_{j}$ of degree less than or equal to $k$, ($k \geq$ 1). \sa
So, we remind the result of \cite{RaTho82} from which our study is developed. Let $\|.\|_{1}$ be the classical norm in $H^1(\Omega)$ and $|.|_{k+1}$ the semi-norm in $H^{k+1}(\Omega)$, and $h$ is the mesh size, namely the largest diameter of the elements of the mesh ${\mathcal T}_h$, then we have:
\begin{lemma}\label{Thm_error_estimate}
Suppose that there exists an integer $k \geq 1$ such that the approximation $u_h$ of $V_h$ is a continuous piecewise function composed by polynomials which belong to $P_k(K_{j}), (1\leq j\leq  N_{s})$. \sa
Then, if the exact solution $u$ to (\ref{VP}) belongs to $H^{k+1}(\Omega)$, we have the following error estimate:
\begin{equation}\label{estimation_error}
\|u_h-u\|_{1} \hs \leq \hs \mathscr{C}_k\,h^k \, |u|_{k+1}\,,
\end{equation}
where $\mathscr{C}_k$ is a positive constant independent of $h$.
\end{lemma}
Let us now consider two families of Lagrange finite elements $P_{k_1}$ and $P_{k_2}$ corresponding to a set of values $(k_1,k_2)\in \N^2$ such that $0< k_1 < k_2$. \\[0.1cm]
The two corresponding inequalities given by (\ref{estimation_error}), assuming that the solution $u$ to (\ref{VP}) belongs to $H^{k_2+1}(\Omega)$, are:
\vspace{-0.2cm}
\begin{eqnarray}
\|u^{(k_1)}_h-u\|_{1,\Omega} \hs & \leq & \hs \mathscr{C}_{k_1} h^{k_1}\, |u|_{k_1+1,\Omega}, \label{Constante_01} \\%[0.1cm]
\|u^{(k_2)}_h\hspace{-0.09cm}-u\|_{1,\Omega} \hs & \leq & \hs \mathscr{C}_{k_2} h^{k_2}\, |u|_{k_2+1,\Omega}\,, \label{Constante_02}
%\vspace{-1cm}
\end{eqnarray}
where $u^{(k_1)}_h$ and $u^{(k_2)}_h$ respectively denotes the $P_{k_1}$ and $P_{k_2}$ Lagrange finite element approximations of $u$.\\[0.2cm]
Now, if one considers a given mesh for the finite element of $P_{k_2}$ which would contains whose of $P_{k_1}$ then, for the particular class of problems where the variational formulation (\ref{VP}) is equivalent to a minimization formulation, (see for example \cite{ChaskaPDE}), one can show that the approximation error of $P_{k_2}$ is always lower than the one of $P_{k_1}$, and $P_{k_2}$ is more accurate than $P_{k_2}$ for all values of the mesh size $h$.\sa
Then, for a given mesh size value of $h$, we consider two independent meshes for $P_{k_1}$ and $P_{k_2}$ built by a mesh generator. So, usually, to compare the relative accuracy between these two finite elements, one asymptotically considers inequalities (\ref{Constante_01}) and (\ref{Constante_02}) to conclude that, when $h$ goes to zero, $P_{k_2}$ finite element is more accurate that $P_{k_1}$, since $h^{k_2}$ goes faster to zero than $h^{k_1}$. \sa
However, for any application $h$ has a fixed value and this way of comparison is not valid anymore. Therefore, our point of view will be to determine the relative accuracy between two finite elements $P_{k_1}$ and $P_{k_2}, (k_1<k_2)$, for any given value of $h$ for which two independent meshes have to be considered.\sa
To this end, let us set:
\begin{equation}\label{beta_ki}
\beta_{k_1} = \mathscr{C}_{k_1}h^{k_1} |u|_{k_1+1,\Omega} \mbox{ and } \beta_{k_2} = \mathscr{C}_{k_2}h^{k_2} |u|_{k_2+1,\Omega}.
\end{equation}
Therefore, instead of (\ref{Constante_01}) and (\ref{Constante_02}), we consider in the sequel the two next inequalities:
\begin{eqnarray}
\|u^{(k_1)}_h-u\|_{1,\Omega} & \leq & \beta_{k_1}, \label{Constante_01_2} \\%[0.1cm]
\|u^{(k_2)}_h\hspace{-0.09cm}-u\|_{1,\Omega} & \leq & \beta_{k_2}. \label{Constante_02_2}
%\vspace{-1cm}
\end{eqnarray}
Now, as we explained in \cite{CMAM1}, there is no {\em a priori} available information to surely or better specify the relative position between the aproximation errors $\|u^{(k_1)}_h-u\|_{1,\Omega}$ and $\|u^{(k_2)}_h-u\|_{1,\Omega}$ which respectively live in the interval $[0, \beta_{k_1}]$ and $[0, \beta_{k_2}]$. \sa
Moreover, we also motivated in \cite{CMAM1} that we have to deal with finite element methods where quantitative uncertainties have to be taken into account in their calculations. This mainly comes from the way the mesh grid generator will process the mesh to compute the approximation $u^{(k_i)}_h, (i=1,2)$, leading to a partial non control of the mesh, even for a given maximum mesh size. As a consequence, the corresponding grid is \emph{a priori} random, and the corresponding approximation $u^{(k_i)}_h, (i=1,2),$ too. \sa
\noindent For these reasons, let us recall the convenient probabilistic framework we introduced in \cite{CMAM1} to consider the possible values of the norm $\|u^{(k)}_h-u\|_{1,\Omega}$ viewed as a random variable defined as follows:
\begin{itemize}
\item For a fixed value of the mesh size $h$, a {\em random trial} corresponds to the grid constitution and the associated approximation $u^{(k)}_h$.
\item The probability space ${\bf\Omega}$ contains therefore all the possible results for a given random trial, namely, all of the possible grids that the mesh generator may processed, or equivalently, all of the corresponding associated approximations $u^{(k)}_h$.
\end{itemize}
Then, for a fixed value of $k$, we define by $X^{(k)}$ the random variable as follows:
\begin{eqnarray}
X^{(k)} : & {\bf\Omega} & \hspace{0.1cm}\rightarrow \hspace{0.2cm}[0,\beta_k] \noindent \\%[0.2cm]
& \boldsymbol{\omega}\equiv u^{(k)}_h & \hspace{0.1cm} \mapsto \hspace{0.2cm}\D X^{(k)}(\boldsymbol{\omega}) = X^{(k)}(u^{(k)}_h) = \|u^{(k)}_h-u\|_{1,\Omega}. \label{Def_Xi_h}
\end{eqnarray}
In the sequel, for simplicity, we will set: $X^{(k)}(u^{(k)}_h)\equiv X^{(k)}(h)$. \sa
So, our interest is to evaluate the probability of the event
\begin{equation}\label{objectif}
\left\{\|u^{(k_2)}_h-u\|_{1,\Omega} \leq \|u^{(k_1)}_h-u\|_{1,\Omega}\right\} \equiv \left\{X^{(k_2)}(h) \leq X^{(k_1)}(h)\right\},
\end{equation}
which will enable us to estimate the more likely accurate between two finite elements of order $k_1$ and $k_2$, $(k_1<k_2)$.\sa
Now, in \cite{ArXiv_JCH}, regarding the absence of information concerning the more likely or less likely values of the norm $\|u^{(k_i)}_h-u\|_{1,\Omega}, (i=1,2),$ in the interval $[0, \beta_{k_i}], (i=1,2)$, we assumed that the two random variables $X^{(k_i)}, (i=1,2),$ have a uniform distribution on their respective interval $[0, \beta_{k_i}]$, and also, that they are independent as well.\sa
Then, this probabilistic framework enabled us to get in a more general context, (see Theorem 3.1 in \cite{ArXiv_JCH}), the density of probability of the random variable $Z$ defined by $Z=X^{(k_2)}-X^{(k_1)}$, and as a consequence, the ${\cal P}_{k_1,k_2} \equiv Prob\left\{ X^{(k_2)} \leq X^{(k_1)}\right\}$ which corresponds to the value of the entire cumulative distribution function $F_Z(z)$ at $z=0$ defined by:
\begin{equation}\label{fonction repartition_FZ}
F_Z(z) = \int_{-\infty}^{z}f_Z(z) dz.
\end{equation}
For the purpose of the present work, the same results may be obtained by elementary adaptations of the results of Theorem 3.3 in \cite{ArXiv_JCH} to get the corresponding probability law given by:
\begin{flushleft}
$\D \hspace{1.5cm}{\cal P}_{k_1,k_2}(h) \, = \,\left| \mbox{\hfill \begin{minipage}[h]{8cm}
\vspace{-0.4cm}
\begin{eqnarray}
\hspace{-0.5cm}\D 1-\frac{1}{2}\!\left(\!\frac{\!\!h}{h^{*}_{k_1,k_2}}\!\right)^{\!\!k_2-k_1} & \hspace{-0.2cm}\mbox{ if }\hspace{-0.2cm} & \hs 0 \leq h \leq h^{*}_{k_1,k_2}, \label{CMAM1}\\[0.2cm]
\hspace{-0.5cm}\D \frac{1}{2}\!\left(\!\frac{h^{*}_{k_1,k_2}}{\!\!h}\!\right)^{\!\!k_2-k_1} & \hspace{-0.2cm}\mbox{ if }\hspace{-0.2cm} & h \geq h^{*}_{k_1,k_2},\label{CMAM2}
\end{eqnarray}
\end{minipage}
}
\right.$
\end{flushleft}
where $h^{*}_{k_1,k_2}$ is defined by:
\begin{equation}\label{h*}
\D h^{*}_{k_1,k_2} \equiv\left( \frac{\mathscr{C}_{k_1}|u|_{k_1+1,\Omega}}{\mathscr{C}_{k_2}|u|_{k_2+1,\Omega}}\right)^{\frac{1}{k_2-k_1}}.
\end{equation}
The shape of this law looks like to a "sigmoid" curve as one can see in Figure \ref{Sigmoid_Curve}.\sa
We already remarked in \cite{CMAM2} that the probability law (\ref{CMAM1})-(\ref{CMAM2}) can asymptotically - when $k_2-k_1$ goes to infinity - leads to the limit situation we defined as the "two-steps" model, (see Figure \ref{Sigmoid_Curve}). \sa
But, we also proved in \cite{CMAM1} that under suitable probabilistic assumptions, one can directly derive this "two-steps" probabilistic law for all non zero integers $k_1$ and $k_2$ to finally get the following probability law:
\begin{equation}\label{Heaviside_Prob}
\D {\cal P}_{k_1,k_2}(h)= \left |
\begin{array}{ll}
\hs 1 & \mbox{ if } \hs 0 < h < h^{*}_{k_1,k_2}, \medskip \\
\hs 0 & \mbox{ if } \hs h> h^{*}_{k_1,k_2}.
\end{array}
\right.
\end{equation}
The next section is dedicated to the analysis of the fit between statistical data and the above probability laws (\ref{CMAM1})-(\ref{CMAM2}) and (\ref{Heaviside_Prob}).
\section{Comparison between numerical statistics and the relative error accuracy probability law}\label{First_Comparison}
\begin{figure}[htb]
  \centering
  \includegraphics[width=12.cm]{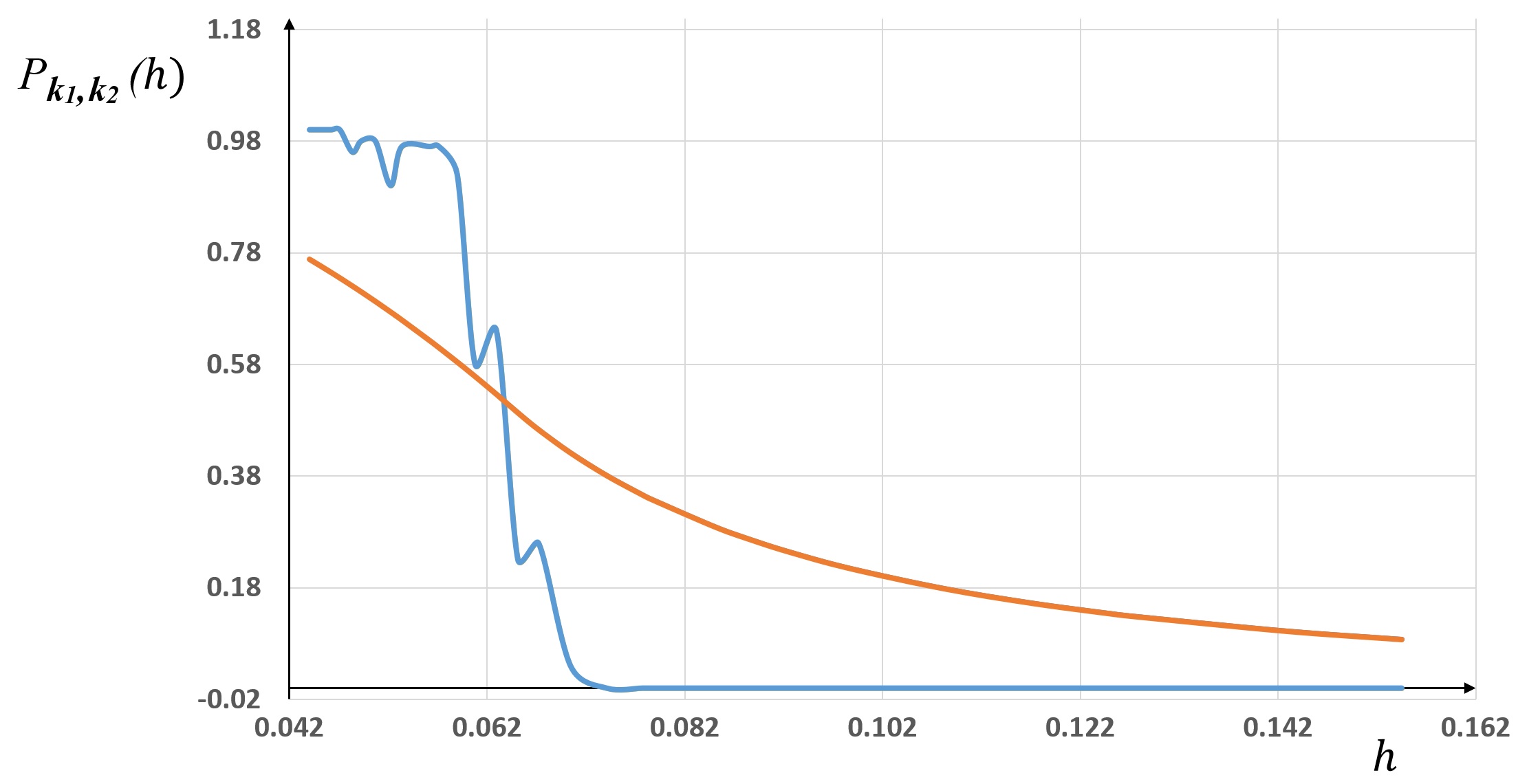}
  \caption{$P_1$ versus $P_3$ for the Runge function with $\alpha=3000$. Comparison between the statistical frequencies (blue) and the "Sigmoid law (orange)} \label{P1P3alpha3000}
\end{figure}

\noindent This section is devoted to the comparison between the theoretical probability law (\ref{CMAM1})-(\ref{CMAM2}) and the statistics one can get from a particular case of the variational formulation (\ref{VP}). \sa
Indeed, since we already showed in \cite{Fitting Stat_Heaviside}, the "two-steps" law (\ref{Heaviside_Prob}) fits well several numerical cases. However this law is a bit rough and cannot really follow the variations of the convexity of the statistical data, (see below). It is the reason why we also tested the accuracy of the fit with the "sigmoid" law given by (\ref{CMAM1})-(\ref{CMAM2}).
More precisely, we considered numerical approximations we implemented for the so-called Poisson-Dirichlet partial differential equation defined in the open unit square $\Omega$ of $\R^2$, where we choose to build the solution by the help of the famous Runge function given by $$f(t)=\D\frac{1}{1+\alpha t^2},$$ where $\alpha$ is a real parameter, (\cite{Rossi}, \cite{Eppe87}). \sa
So, to numerically check the accuracy of the fit between a given probability law and the corresponding statistical data produced by numerical simulations, we considered two finite element $P_{k_1}$ and $P_{k_2}$ ($k_1<k_2$) and we fixed a given number of meshes to be built by the generator of meshes, each of them associated to a same given mesh size $h$. \sa
Then, to evaluate the relative accuracy between the two concerned finite elements, we tested for each mesh if $\|u^{(k_2)}_h-u\|_{1}$ is lower than $\|u^{(k_1)}_h-u\|_{1}$. Then, we repeated the same process for different values of $h$ which gave us a function of $h$, namely, the frequency of cases when the approximation error associated to $P_{k_2}$ finite element is lower than the one computed with the $P_{k_1}$ one. In all cases, we use the FreeFem++ package~\cite{Hech12} to compute the $P_{k_i}, (i=1,2),$ finite element approximations. \sa
\noindent Now, to motivate the next section, let us recall a typical result we got in \cite{Fitting Stat_Heaviside}. To this end, we consider here the particular numerical test we implemented to analyze the relative accuracy between the $P_1$ and $P_3$ finite elements. We considered 100 meshes for each value of $h$ and we fixed the parameter $\alpha$ of the Runge function to value of 3000. \sa
Then, in Figure \ref{P1P3alpha3000}, for the different concerned values of the mesh size $h$, let us plot on the results obtained by the statistical frequencies corresponding to the cases such that the finite element $P_{k_2}$ is more accurate that the $P_{k_1}$ one (red), together with the corresponding "sigmoid" probability law (blue).\sa
The "sigmoid" curve showed in Figure \ref{P1P3alpha3000} was computed by the help of the Excel solver to statistically determine the value of $h^{*}_{k_1,k_2}$ by a least squares adjustment between the two concerned curves.\sa
As one can see, the fit between the statistics and the probability law is not satisfactory. The same gaps were observed for other simulations corresponding to different sets of parameter values and for several pairs of finite elements $P_{k_1}$ and $P_{k_2}$. \sa
So, one could expect to get a better fit that the one provided by the "Sigmoid" probability law (\ref{CMAM1})-(\ref{CMAM2}). This is the purpose of the next section where we will show how to enrich the "sigmoid" which only depends on one parameter, namely, $h^{*}_{k_1,k_2}$.
\section{The new probability law for the relative error accuracy between two finite elements $P_{k_1}$ and $P_{k_2}, (k_1<k_2).$}\label{New probabilistic_law}
\noindent This section is devoted to the new probabilistic law we will derive to evaluate the relative error accuracy between two finite elements $P_{k_1}$ and $P_{k_2}, (k_1<k_2)$.
To motivate the new position we will consider, let us proceed to a couple of remarks:
\begin{enumerate}
\item The first one concerns the assumptions we took into account in the previous works that enabled us to derive the "Sigmoid" probability law. Since we would like to get a more precise fit between the probabilistic law and the statistical data, we will relax the hypothesis of uniformity we applied to the densities of the random variables $X^{(k_i)}(h), (i=1,2)$, (see \cite{ArXiv_JCH}, \cite{CMAM1} and \cite{ChAs20}).
\item To choose the shape of these densities, first of all, we will consider the one of the random variable $Z=X^{(k_2)}-X^{(k_1)}$. Indeed, since our goal is to get for the cumulative distribution function $F_Z$ defined by (\ref{fonction repartition_FZ}), at the point $z=0$, a curve whose shape looks like a "Sigmoid", we will enrich our modeling process by adding more degrees of freedom. More precisely, if the "Sigmoid" probability law (\ref{CMAM1})-(\ref{CMAM2}) contains one parameter which is $h^{*}_{k_1,k_2}$, we will now consider a density $f_Z$ for the random variable $Z$ such that the corresponding value of its cumulative distribution function $F_Z$ at the point $z=0$ will include two exogenous parameters to be statistically estimated.
\end{enumerate}
Keeping in mind these remarks, we begin by introducing the probability density function $f_X$ of the normalized Beta random variable $X$ defined by:
%!
\begin{equation}
\D\forall x \in [0,1]: f_X(x;p,q) \equiv \frac{x^{p-1}(1-x)^{q-1}}{\D\int_{0}^{1}u^{p-1}(1-u)^{q-1}du} = \frac{\Gamma(p+q)}{\Gamma(p)\Gamma(q)}x^{p-1}(1-x)^{q-1},
\end{equation}
where $p$ and $q$ are two parameters of shape which belong to $\R^{*}_{+}$ and $\Gamma(.)$ denotes the classical Gamma function. \sa
Among the numerous features of the Beta distribution, let us mention one of them which particularly motivates us to consider it for our objectives. Indeed, depending on the two parameters $p$ and $q$, the shapes of its corresponding cumulative distribution are very rich and include the shape of the "Sigmoid" curve we are looking for fitting the statistics when one considers the case of the Runge solution to the Laplacian-Dirichlet problem in the unit square of $\R^2$.\sa
However, one cannot directly apply the Beta distribution to get the probability law we are looking for. More precisely, two main features have to be taken into account:
\begin{enumerate}
\item If the support of the Beta density $f_X$ of the random variable $X$, denoted $Supp X$, is include in $[0,1]$, the one of the random variable $Z$ is $[-\beta_{k_1}, \beta{k_2}]$, since $Z=X^{(k_2)}-X^{(k_1)}$ and $Supp\,X^{(k_i)}\subset [0,\beta_{k_i}], (i=1;2)$. This will drive us to a suitable transformation of the density $f_X$ to guarantee the correct support of the density $f_Z$ of $Z$.
\item Given that we are looking for a probability law of the event $\D\left\{X^{(k_2)} \leq X^{(k_1)}\right\}$ as a function of $h$ which belongs to $[0, +\infty[$, then we will also apply consequently another transformation of the density $f_Z$ to assure this property for the support of $h$.
\end{enumerate}
So, to achieve these transformations we establish the following results:
\begin{lemma}\label{New_f_Z}
Let $Z$ be the random variable defined by $Z\equiv X^{(k_2)}-X^{(k_1)}$, where $X^{(k_i)}, (i=1,2),$ are defined by (\ref{Def_Xi_h}). Let $X \sim B(p,q)$ be the Beta distribution parameterized by two given shape parameters $(p,q) \in \R^{*2}$. Assume that $Z$ is defined by:
\begin{equation}\label{Decomp_Z}
Z=-\beta_{k_1}+(\beta_{k_1}+\beta_{k_2})X.
\end{equation}
Then, the probability density function $f_Z$ of the $Z$ is given by:
\begin{equation}\label{f_Z}
\D f_{Z}(z)= \frac{\Gamma(p+q)}{\Gamma(p)\Gamma(q)}\,\frac{\beta_{k_1}^{p-1}\beta_{k_2}^{q-1}}{(\beta_{k_1}+\beta_{k_2})^{p+q-1}}\left(1+\frac{z}{\beta_{k_1}}\right)^{\!p-1}
\!\!\left(1-\frac{z}{\beta_{k_2}}\right)^{\!q-1}\mathbbm{1}_{[-\beta_{k_1},\beta_{k_2}]}(z),
\end{equation}
where $\mathbbm{1}_{[-\beta_{k_1},\beta_{k_2}]}$ is the indicator function of the interval $[-\beta_{k_1},\beta_{k_2}]$.
\end{lemma}
\begin{prooff}
Since we already noticed, the support of $f_{Z}$ is clearly within $[-\beta_{k_1},\beta_{k_2}]$ as soon as those of $X^{(k_i)}$ are in $[0,\beta_{k_i}]$.\sa
Now, let us evaluate the cumulative distribution function $F_Z(z)$ of the random variable $Z$ defined by (\ref{Decomp_Z}):
\begin{eqnarray}
\D F_Z(z) & = & \D Prob\left\{\frac{}{}\!Z\leq z\right\} = Prob\left\{\frac{}{}\!\!-\beta_{k_1}+(\beta_{k_1}+\beta_{k_2})X \leq z\right\} \\[0,2cm]
 & = & Prob\left\{\frac{}{}\!\!X \leq \D\frac{z+\beta_{k_1}}{\beta_{k_1}+\beta_{k_2}}\right\} =\frac{\Gamma(p+q)}{\Gamma(p)\Gamma(q)}\int_{0}^{\frac{z+\beta_{k_1}}{\beta_{k_1}+\beta_{k_2}}}u^{p-1}(1-u)^{q-1}du. \label{FZ_to_derive}
\end{eqnarray}
Then, we derive $F_Z(z)$ given by (\ref{FZ_to_derive}) which leads to the expression (\ref{f_Z}) of the density $f_Z$.
\end{prooff}
\begin{remark}
Let us give the meaning of this result. Since we already mentioned, in our previous works (see for example \cite{CMAM1} and \cite{ArXiv_JCH}) we assumed the random variables $X^{(k_i)}$ to be uniformly distributed on their support $[0,\beta_{k_i}]$. This lead us to consequently derived the density $f_Z$ of the random variable $Z$, (see Theorem 3.1 in \cite{ArXiv_JCH}). \sa
Here, Lemma \ref{New_f_Z} may be interpreted about the new assumption we implicitly made on the random variables $X^{(k_i)}$. Indeed, if we rewrite $Z$ like: $$Z=\beta_{k_2}X - (\beta_{k_1}-\beta_{k_1}X),$$
then by setting:
\begin{equation}\label{New_Xki}
X^{(k_1)} \equiv \beta_{k_1} - \beta_{k_1}X \hs \mbox{ and } \hs X^{(k_2)} = \beta_{k_2}X,
\end{equation}
we observe that each support of $X^{(k_i)}$ belongs to $[0,\beta_{k_i}]$, on the one hand, and the difference $X^{(k_2)}-X^{(k_1)}$ is equal to $Z$, on the other hand. \sa
In other words, the choice we consider here to write the random variable $Z$ given by (\ref{Decomp_Z}) as a dimensional Beta distribution on $[-\beta_{k_1}, \beta_{k_2}]$, (which corresponds to alter the location and scale of the standard Beta distribution), consequently modifies the hypothesis of uniformity of the two variables $X^{(k_i)}$ to the one of a non dimensional Beta distributions.
\end{remark}
We are now in position to derive the new probability distribution ${\cal P}_{k_1,k_2}(h)$ to evaluate the more likely accurate between two Lagrange finite elements $P_{k_1}$ and $P_{k_2}, (k_1<k_2)$.
This is the purpose of the following Theorem.
\begin{theorem}\label{New_GBP_Law_Th}
Let $X^{(k_i)}, (i=1,2),$ be the two random non dimensional Beta distribution variables defined by (\ref{New_Xki}) and $Z$ the corresponding random variable defined by (\ref{Decomp_Z}). \sa
Then, ${\cal P}_{k_1,k_2}(h)$ is the cumulative distribution function of a generalized Beta prime random variable $H$, whose density of probability $f_H$, is defined by four parameters $(p,q,k_2-k_1, h^{*}_{k_1,k_2})$, and we have:
\begin{equation}\label{New_Law}
\D{\cal P}_{k_1,k_2}(h) = Prob\{H\geq h\} = \int_{h}^{+\infty} f_H(s;q,p,k_2-k_1, h^{*}_{k_1,k_2})\, ds,
\end{equation}
where:
\begin{equation}\label{fH}
\D f_H(s;q,p,k_2-k_1, h^{*}_{k_1,k_2}) = \frac{\Gamma(p+q)}{\Gamma(p)\Gamma(q)} \, \frac{(k_2-k_1)}{h^{*}_{k_1,k_2}} \, \left(\frac{s}{h^{*}_{k_1,k_2}}\!\right)^{q(k_2-k_1)-1}\left[1+\left(\frac{s}{h^{*}_{k_1,k_2}}\!\right)^{k_2-k_1}\right]^{-p-q}.
\end{equation}
\end{theorem}
\begin{prooff}
Since we are looking for a random variable $H$ whose support has to be $[0,+\infty[$, which is associated to ${\cal P}_{k_1,k_2}(h)$ equals to $F_Z(0)$, on the one hand, and (\ref{FZ_to_derive}) was derived by the help of the dimensionless Beta distribution $B(p,q)$ whose support belongs to $[0,1]$, on the other hand, we set the following change of variable in (\ref{FZ_to_derive}):
$$ \D s=\frac{u}{1-u},$$
and we get:
\begin{equation}\label{F-Z_to_derive_2}
\D F_{Z}(z) =\frac{\Gamma(p+q)}{\Gamma(p)\Gamma(q)}\int_{0}^{\frac{z+\beta_{k_1}}{\beta_{k_2}-z}}s^{p-1}(1+s)^{-p-q}\,ds.
\end{equation}
Now, let us set in (\ref{F-Z_to_derive_2}) $z=0$. Then, by using $\beta_{k_i}, (i=1,2),$ given by (\ref{beta_ki}) and $h^{*}_{k_1,k_2}$ by (\ref{h*}), we obtain:
\begin{eqnarray}
\D F_{Z}(0) \hs = \hs {\cal P}_{k_1,k_2}(h) & = & \D \frac{\Gamma(p+q)}{\Gamma(p)\Gamma(q)}\int_{0}^{\frac{\beta_{k_1}}{\beta_{k_2}}}s^{p-1}(1+s)^{-p-q}\,ds, \\[0.2cm]
& = & \frac{\Gamma(p+q)}{\Gamma(p)\Gamma(q)}\int_{0}^{(h^{*}_{k_1,k_2}/h)^{k_2-k_1}}s^{p-1}(1+s)^{-p-q}\,ds.\label{F-Z_to_derive_3}
\end{eqnarray}
A last change of variables $t=\frac{1}{s}$ in (\ref{F-Z_to_derive_3}) leads to:
\begin{equation}\label{F-Z_to_derive_4}
\D F_{Z}(0) = \frac{\Gamma(p+q)}{\Gamma(p)\Gamma(q)} \int_{(h/h^{*}_{k_1,k_2})^{k_2-k_1}}^{+\infty} t^{q-1}(1+t)^{-p-q}\,dt.
\end{equation}
Finally, by considering  the probability of the complementary event considered in (\ref{F-Z_to_derive_4}) as a function of $h$, (namely, $1-{\cal P}_{k_1,k_2}(h)$), after derivation with respect to $h$, we get the density $f_H$ defined in (\ref{fH}) which is a generalized prime Beta probability density function parameterized by $(q,p,k_2-k_1, h^{*}_{k_1,k_2})$, (see for example \cite{Gavin_E_Crooks}). \sa
In other words, we have:
\begin{equation}\label{Equiv}
Prob\left\{X^{(k_1)}\leq X^{(k_2)}\right\} = 1 - Prob\left\{X^{(k_2)}\leq X^{(k_1)}\right\}(h) = \int_{o}^{h}f_H(s;q,p,k_2-k_1, h^{*}_{k_1,k_2}) ds,
\end{equation}
and $\D{\cal P}_{k_1,k_2}(h)$ in (\ref{New_Law}) is deduced by complementarity.
\end{prooff}
\begin{remark}$\frac{}{}$
\begin{enumerate}
\item From (\ref{New_Law})-(\ref{fH}), or equivalently by (\ref{F-Z_to_derive_3}), we observe the asymptotic behavior of the probability of the event $\D\left\{X^{(k_1)}\leq X^{(k_2)}\right\}$ when $h$ when the mesh size $h$ goes to 0. \sa
    Clearly, it goes to 0 when $h$ goes itself to 0. In other words, we found with the new probabilistic law the classical result which claims that $P_{k_2}$ finite element is always more accurate the $P_{k_1}$ one, since from (\ref{Constante_01}) and (\ref{Constante_02}), $h^{k_2}$ goes faster to zero than $h^{k_1}$ when $k_1<k_2$. \sa
    Here, by (\ref{F-Z_to_derive_3}) the same property is expressed in terms of probability, namely: the event $\D\left\{X^{(k_1)}\leq X^{(k_2)}\right\}$ is an almost never one, or equivalently, the event $\D\left\{X^{(k_2)}\leq X^{(k_1)}\right\}$ is an almost surely one since its probability is equal to 1.
    \item From (\ref{New_Law}) and since the positivity of the density $f_H$, we conclude that $\D{\cal P}_{k_1,k_2}(h)$ is a decreasing function of $h$. This property was already observed with the "Sigmoid" model in \cite{CMAM1} and \cite{ArXiv_JCH} when uniformity of the random variables $X^{(k_i)}, (i=1,2),$ was assumed. \sa
        In other words, this property claims that the event $\D\left\{P_{k_2} \mbox{ is more accurate than } P_{k_1}\right\}$ is not an absolute reality and the more $h$ increases the less this property is less likely. Moreover, when $h$ becomes great, asymptotically the event $\D\left\{P_{k_2} \mbox{ is more accurate than } P_{k_1}\right\}$ in an almost never one !
\end{enumerate}
\end{remark}
\section{Numerical comparison between the Generalized Beta prime probabilistic law and statistical frequencies}\label{New_Stat}\label{Stat_Proba}
\begin{figure}[htbp!]
\begin{tabular}{lr}
{
\hspace*{-0.5cm}
\includegraphics[width=8cm]{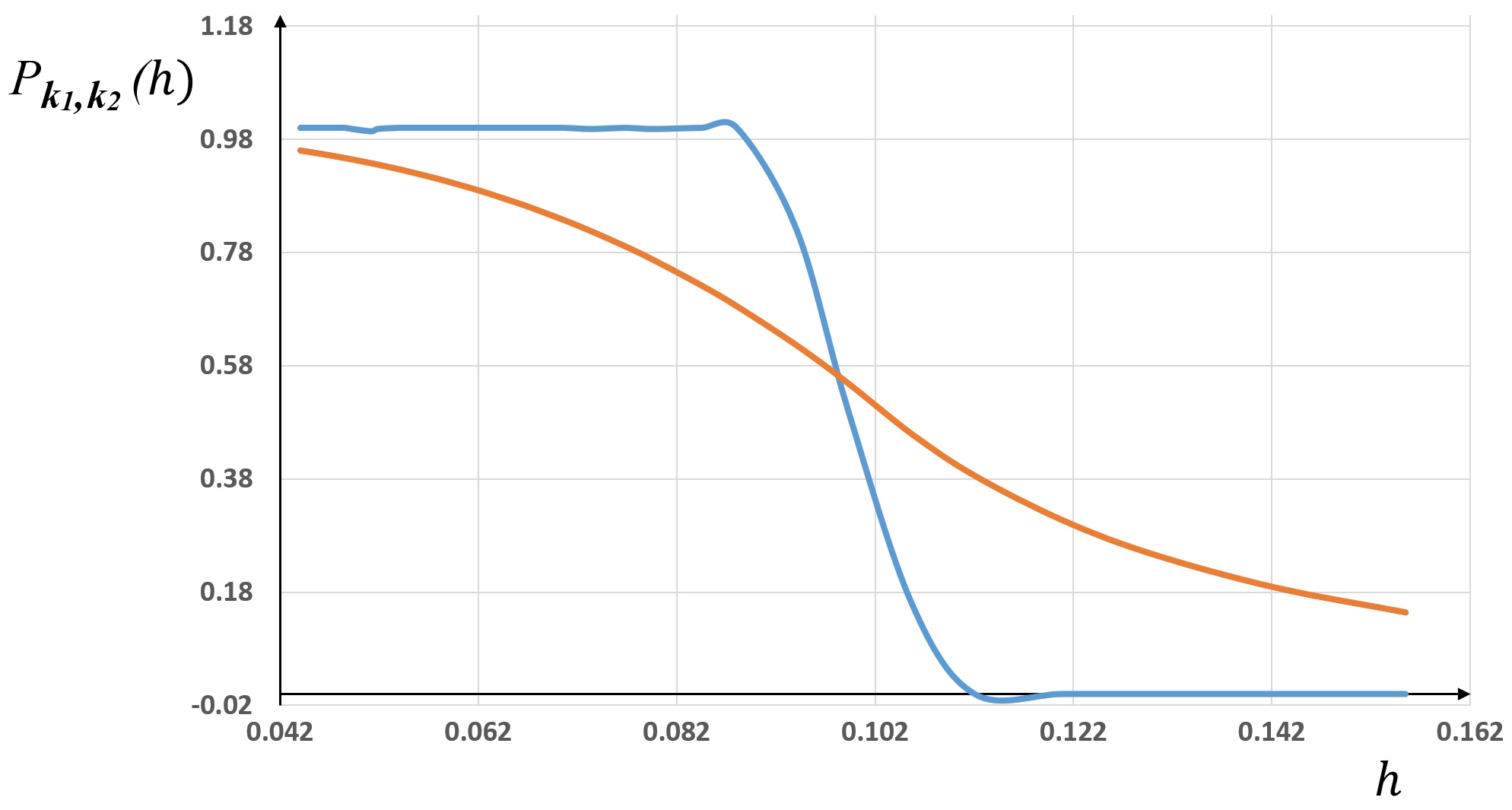}
}
&
\hspace*{-0.5cm}
{
\includegraphics[width=8cm]{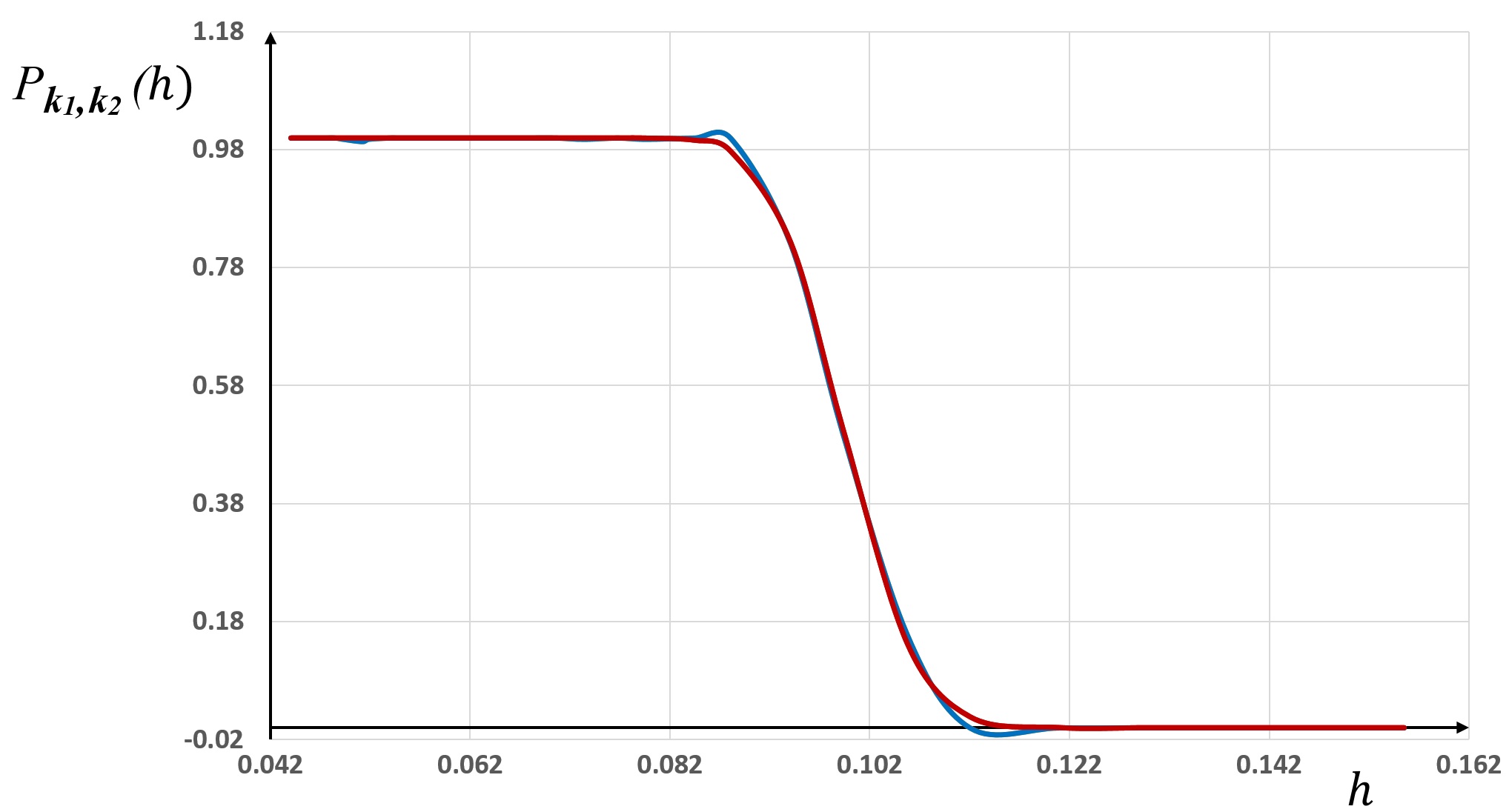}
}
\end{tabular}
\caption{$P_1$ versus $P_4$, Sigmoid (left) and Generalized Beta Prime (right).} \label{P1P4-Sigmoid vs GBP}
\end{figure}

\begin{figure}[htbp!]
\begin{tabular}{lr}
{
\hspace*{-0.5cm}
\includegraphics[width=8cm]{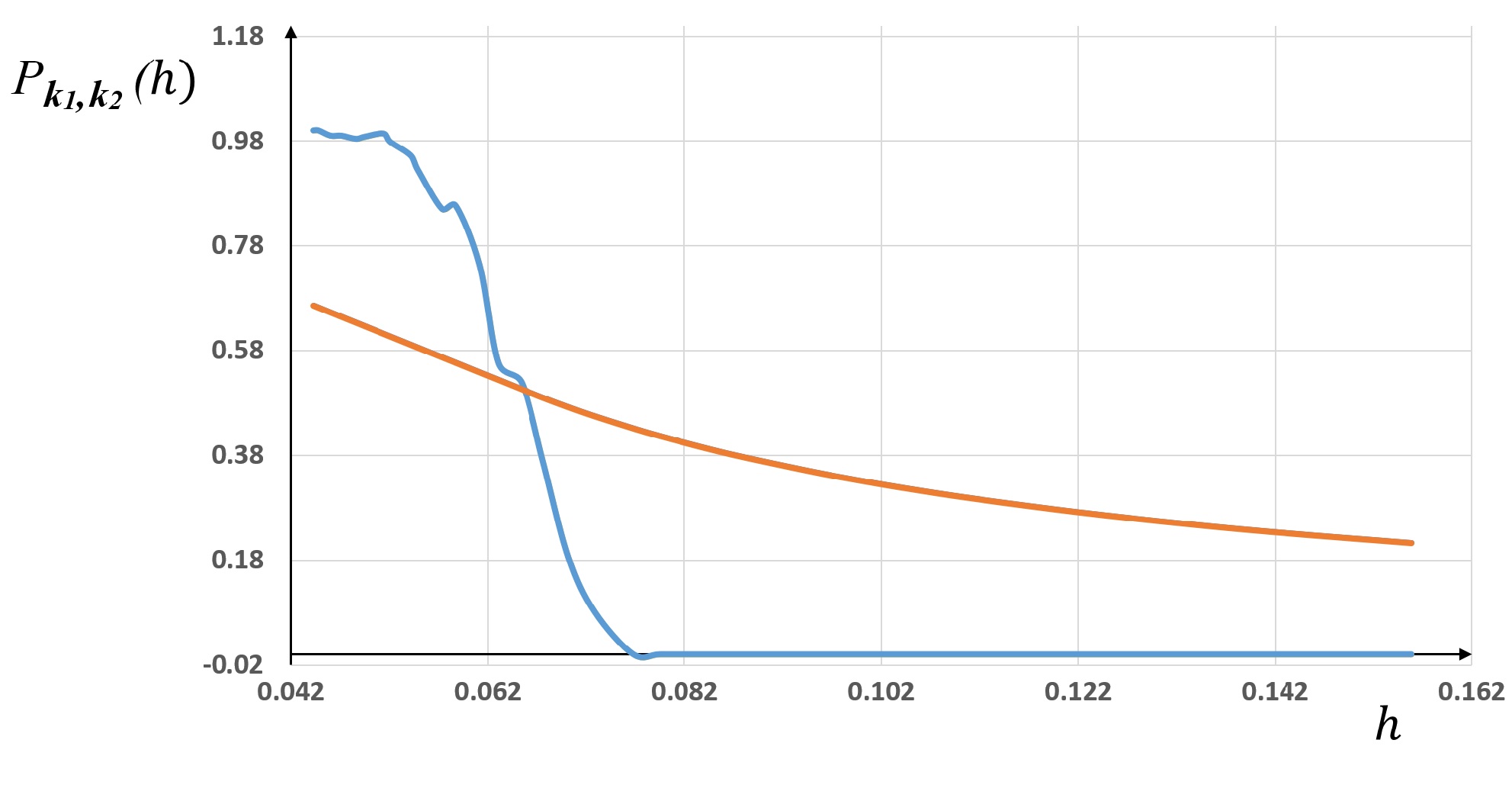}
}
&
\hspace*{-0.5cm}
{
\includegraphics[width=8cm]{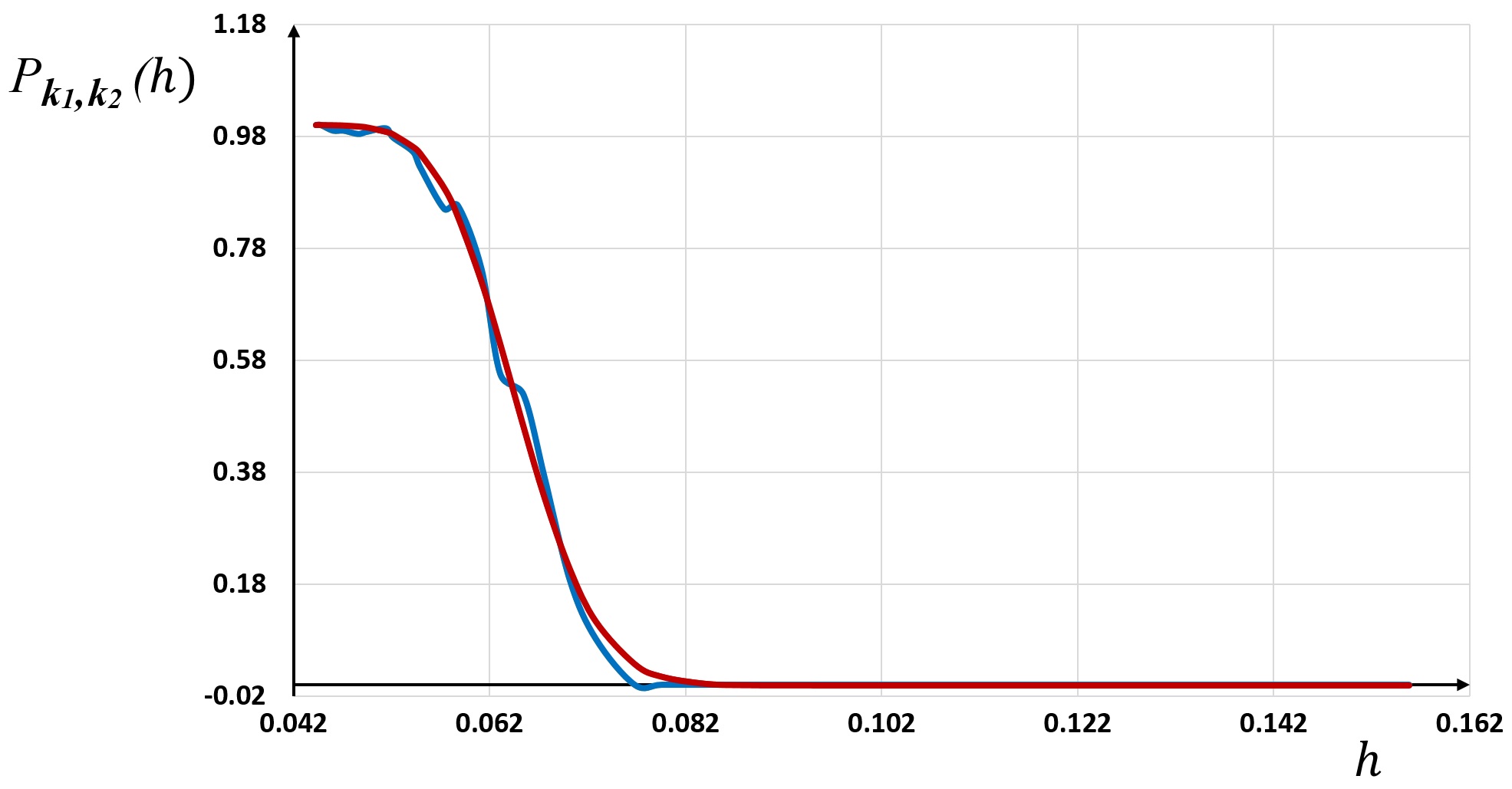}
}
\end{tabular}
\caption{$P_2$ versus $P_3$, Sigmoid (left) and Generalized Beta Prime (right).} \label{P2P3-Sigmoid vs GBP}
\end{figure}

\noindent This section in devoted to evaluate the quality of the fit between the statistical frequencies we presented in Section \ref{First_Comparison} processed by the help of the Runge function and the corresponding probabilities computed by the help of the generalized Beta prime distribution (\ref{New_Law}) that we derived in Theorem \ref{New_GBP_Law_Th}.\sa
To this end, let us explain how we processed to determine the four parameters $(q,p,k_2-k_1,h^{*}_{k_1,k_2})$ of the density $f_H$. First of all, we considered and fixed two finite elements, for example $P_1$ and $P_3$, and so, $k_2-k_1=2$. \sa
Regarding the three other parameters, $p,q$ and $h^{*}_{k_1,k_2}$, we have implemented again an optimization based on least squares adjustment to minimize the sum of the squares of the differences between the statistical frequencies and the corresponding values computed by the generalized Beta prime model presented above. Once more, this was computed by the help of the Excel solver.\sa
\begin{figure}[htb]
  \centering
  \includegraphics[width=12.cm]{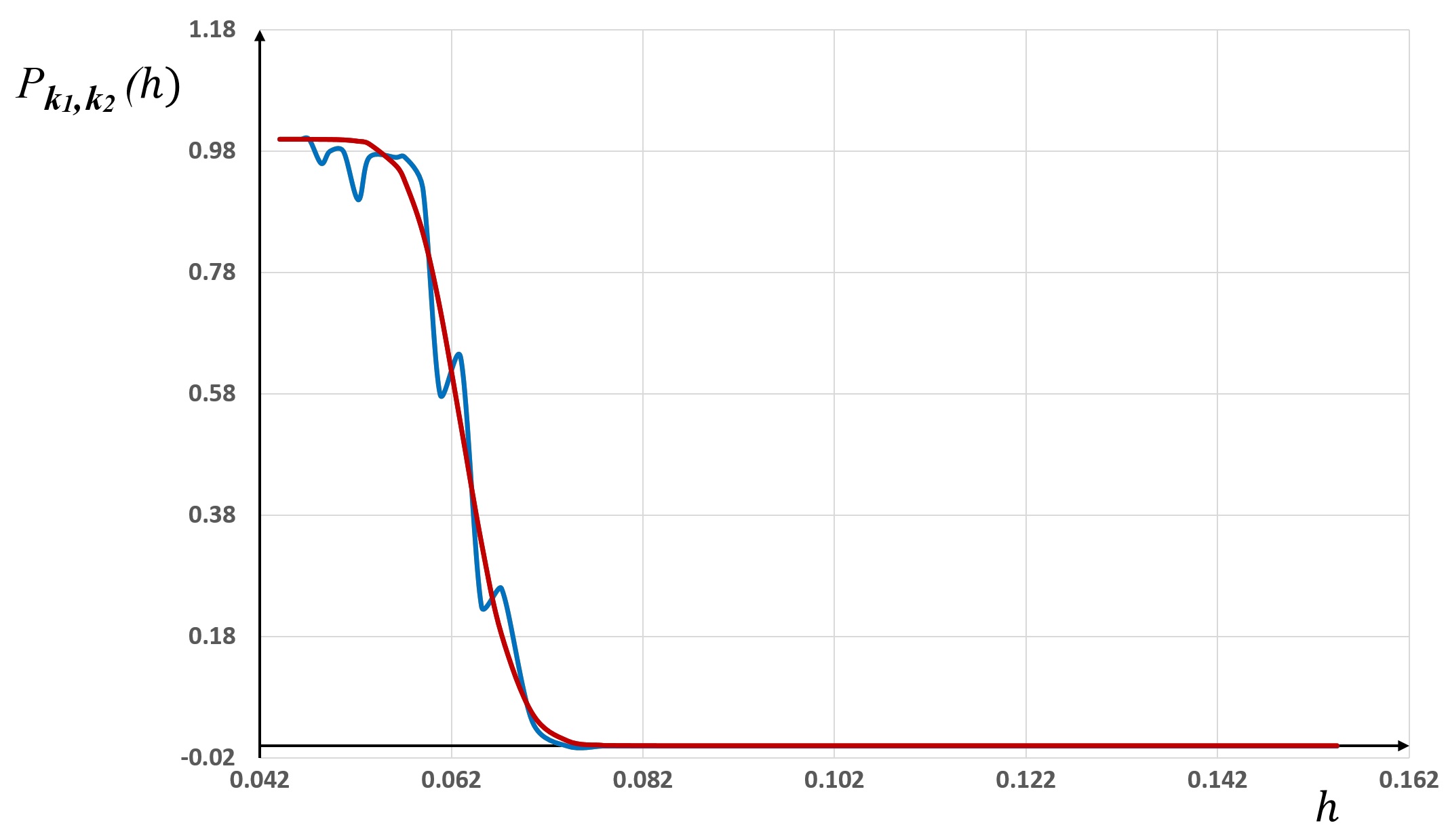}
  \caption{$P_1$ versus $P_3$ for the Runge function with $\alpha=3000$. Comparison between the statistical frequencies (blue) and the Generalized Beta Prime law (red)} \label{Good_Fit}
\end{figure}
Then, we can observe that the least squares algorithm found the optimal parameters $p$, $q$ and $h^{*}_{k_1,k_2}$ such that the fit between the statistical frequencies and the Generalized Beta Prime law (\ref{New_Law}) is very satisfactory, (see Figure \ref{Good_Fit}). \sa
This a clearly due to the the richness and the flexibility of this distribution, providing two degrees of freedom $p$ and $q$ in addition to the "Sigmoid" law, we motivated in Section \ref{New probabilistic_law} for describing the randomness values of the approximation errors $\|u^{(k_i)}_h-u\|_{1,\Omega}, (i=1,2),$ within their respective interval $[0,\beta_{k_i}], (i=1,2)$.\sa
Another example to appreciate the accuracy of the fit may be achieved by comparing the results implemented with the two Lagrange finite elements $P_1$ and $P_4$. In this case, as one can see again in Figure \ref{P1P4-Sigmoid vs GBP}, when the "Sigmoid" law only describes the global trend of the statistical frequencies, the generalized Beta prime law perfectly fits with the corresponding data.\sa
This difference of the quality of the fit between the "Sigmoid" law and the generalized Beta prime law becomes all the more important when one considers $P_2$ finite elements versus $P_3$ one. Indeed, in this case, $k_2-k_1$ is equal to one and the "Sigmoid" law becomes a linear function of $h$ when $h\leq h^{*}_{k_1,k_2}$ in formula (\ref{CMAM1}). On the other hand, one more time, as one can see in Figure \ref{P2P3-Sigmoid vs GBP}, the generalized Beta prime law fit very well with the statistical frequencies.\sa
\section{Conclusions}
\noindent In this paper, we derived a new family of probabilistic laws to compare the accuracy between two Lagrange finite elements $P_{k_1}$ and $P_{k_2}, (k_1<k_2)$. Based on a new probabilistic approach we already introduced in \cite{ArXiv_JCH}, we extended these previous results to improve the fit between the statistical frequencies obtained by implementing practical cases and the corresponding probabilities.\sa
We recall in Section \ref{First_Comparison} the main problems we got regarding the accuracy of the previous probabilistic laws we derived in \cite{ArXiv_JCH},\cite{CMAM1} or \cite{arXiv_Wmp}, and we identify the gaps we observed with the comparable statistics. \sa
Then, in Section \ref{New probabilistic_law} we motivated and derived the new probabilistic law based on the generalized Beta Prime law. Then, we analyzed in Section \ref{Stat_Proba} the quality of the fit we got between the corresponding probabilities which corrected the main observed deficiencies described above.\sa
Finally, this new probabilistic law together with the statistical validation we processed, significantly confirms the relevance to consider the approximation errors like random variables defined in an adapted probabilistic framework. Of course, one must keep in mind that this approach is not limited to finite elements error estimates, but might be fruitful for any kind of error estimates one must deal with other types of numerical approximations. \sa
\textbf{\underline{Homages}:} The authors want to warmly dedicate this research to pay homage to the memory of Professors Andr\'e Avez and G\'erard Tronel who largely promote the passion of research and teaching in mathematics.

\end{document}